\documentclass[12pt]{amsart}
\usepackage{amsmath,amssymb,epsfig,color, amsthm, mathrsfs}
\usepackage{wrapfig,lipsum,floatrow,sidecap,comment}
\usepackage{hyperref}
\usepackage[pagewise]{lineno}
\usepackage{graphicx}
\usepackage{multirow}
\usepackage{hhline}
\usepackage{url}
\graphicspath{ {images/} }
\allowdisplaybreaks

\newtheorem{theorem}{Theorem}[section]

\newtheorem{corollary}[theorem]{Corollary}

\newtheorem{proposition}[theorem]{Proposition}

\newcommand{\vanish}[1]{}\parskip=12pt

\def\p{\prime}

\def\b{\textbf{b}}

\def\K{\mathbb{K}}

\numberwithin{equation}{section}

\newcommand{\bc}{\begin{center}}
\newcommand{\ec}{\end{center}}

\newcommand{\be}{\begin{equation}}
\newcommand{\ee}{\end{equation}}

\newcommand{\beqn}{\begin{eqnarray*}}
\newcommand{\eeqn}{\end{eqnarray*}}
\newcommand{\ev}{\mathrm{ev}}

\theoremstyle{definition}
\newtheorem{notation}{Notation}
\newtheorem{remark}{Remark}

\usepackage{pinlabel}
\newcommand{\pic}[2]{\raisebox{-.5\height}{\includegraphics[scale=#2]{#1}}}

\newcommand\Xor{\pic{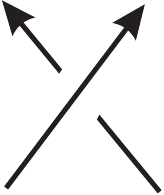}{.50}}
\newcommand\Yor{\pic{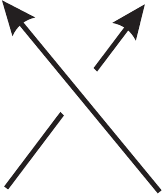} {.50}}
\newcommand\Ior{\pic{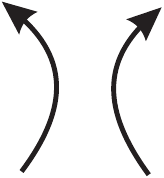} {.50}}
\def\Idor{\pic{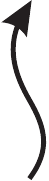} {.50}}

\newcommand\unknot{\pic{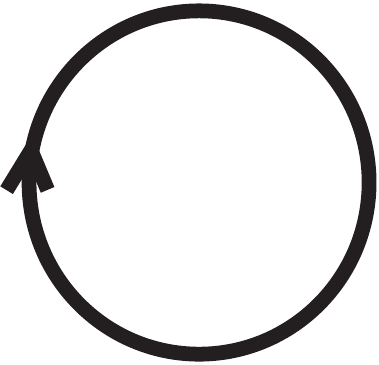} {.20}}

\def\torusparallel{\pic{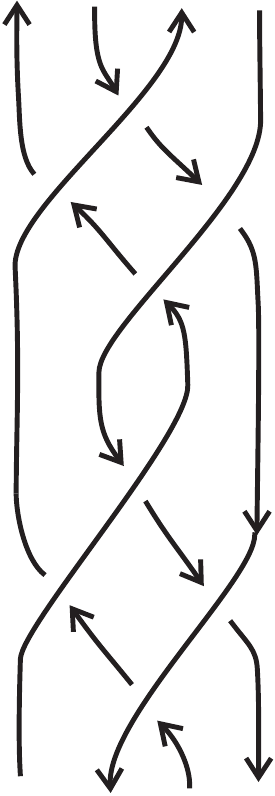} {.50}}
\def\twotwotangle{\pic{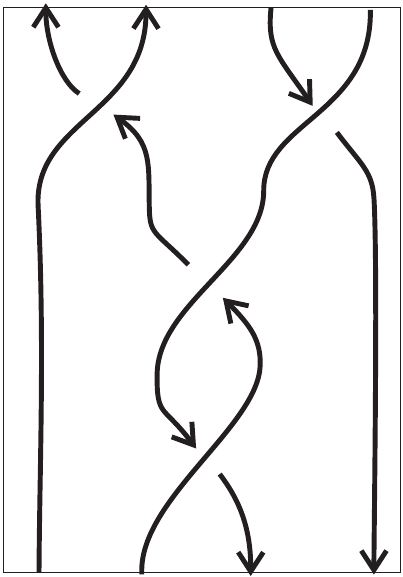} {.4}}
\def\twotwoh{\pic{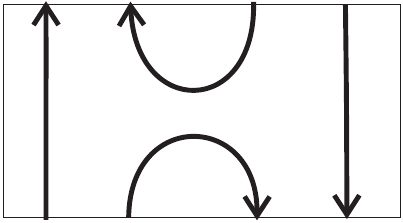} {.4}}
\def\twotwohh{\pic{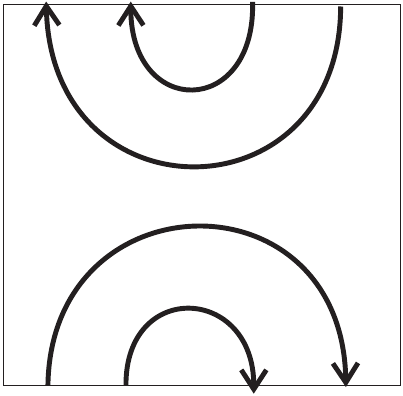} {.4}}
\def\twotwosigmaone{\pic{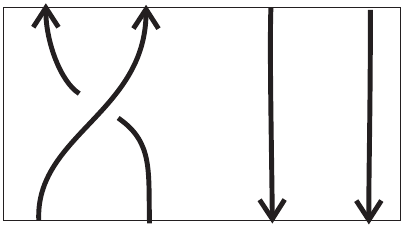} {.4}}
\def\twotwosigmathree{\pic{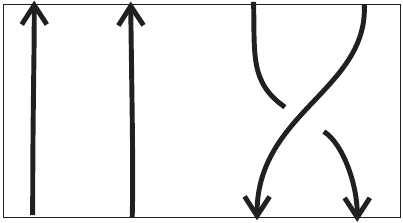} {.4}}
\def\twotwotanglecircle{\pic{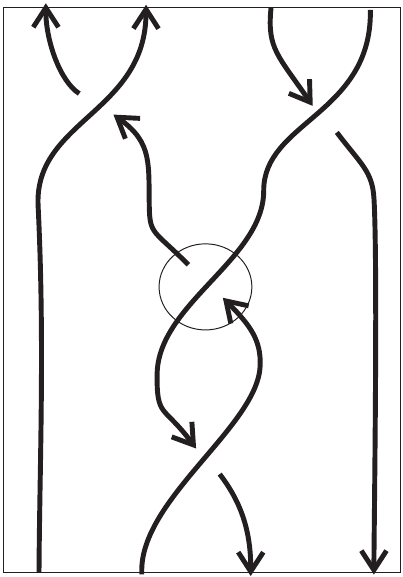} {.4}}
\def\Tswitch{\pic{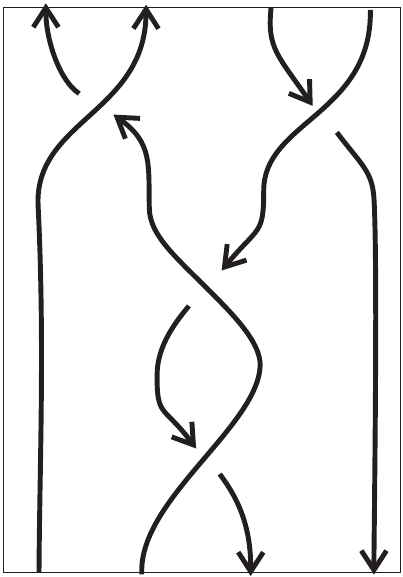} {.4}}
\def\Tsmooth{\pic{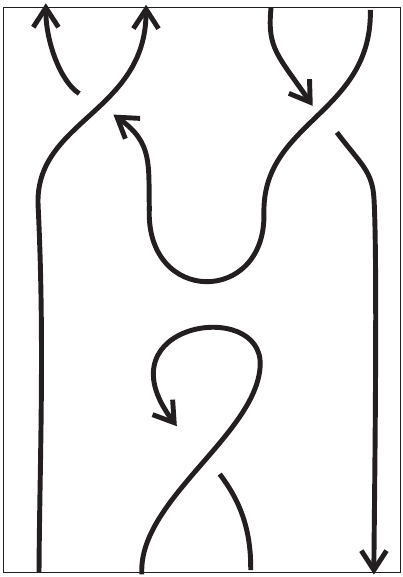} {.4}}
\def\kparallel{\pic{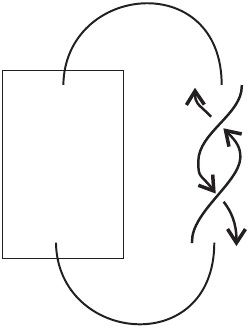}{.7}}
\def\kparallelstraight{\pic{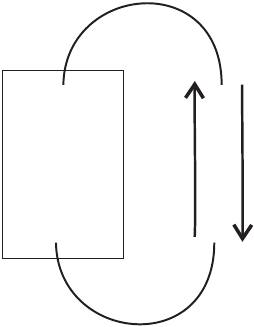}{.7}}
\def\kparallelswitch{\pic{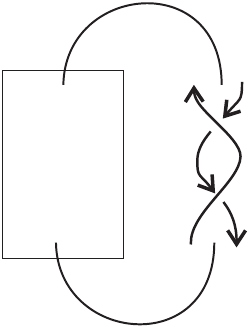}{.7}}
\def\kparallelsmooth{\pic{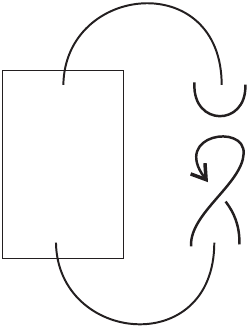}{.7}}
\begin{document}
\title{The Braid Indices of the Reverse Parallel Links of Alternating Knots}
\author{Yuanan Diao$^\sharp$ and Hugh Morton$^\flat$}
\address{$^\sharp$ Department of Mathematics and Statistics\\
University of North Carolina Charlotte\\
Charlotte, NC 28223, USA\\
$^\flat$ Department of Mathematical Sciences\\
University of Liverpool\\
Peach Street, Liverpool. L69 7ZL, United Kingdom }
\email{ydiao@uncc.edu, morton@liverpool.ac.uk}

\begin{abstract}
The braid indices of most links remain unknown as there is no known universal method that can be used to determine the braid index of an arbitrary knot. This is also the case for alternating knots. In this paper, we show that if $K$ is an alternating knot, then the braid index of any reverse parallel link of $K$ can be precisely determined. More precisely, if $D$ is a reduced diagram of $K$, $v_+(D)$ ($v_-(D)$) is  the number of regions in the checkerboard shading of $D$ for which all crossings are positive (negative), $w(D)$ is the writhe of $D$, then the braid index of a reverse parallel link of $K$ with framing $f$, denoted by $\mathbb{K}_f$, is given by the following precise formula
$$
\textbf{b}(\mathbb{K}_f)=\left\{
\begin{array}{ll}
c(D)+2+a(D)-f, &\ {\rm if}\ f < a(D),\\
c(D)+2, &\ {\rm if}\ a(D)\le f \le b(D),\\
c(D)+2-b(D)+f, &\ {\rm if}\ f > b(D),\\
\end{array}
\right.
$$
where $a(D)=-v_-(D)+w(D)$ and $b(D)=v_+(D)+w(D)$. 
\end{abstract} 

\maketitle

\section{Introduction}\label{s1} The determination of the braid index of a knot or a link is known to be a challenging problem. To date there is no known method that can be used to determine the precise braid index of an arbitrarily given knot/link. This is also the case when we restrict ourselves to alternating knots and links, although the braid indices of many alternating knots and links can now be determined, for example all two bridge links and all alternating Montesinos links \cite{Diao2021,Murasugi91}. However, we shall prove in this paper a somewhat surprising result, that is, the braid index of any {\em reverse parallel} link of an alternating knot can be precisely determined. Furthermore, the formula can be derived easily from any reduced diagram of the alternating knot.

In this paper we shall study the {\em reverse parallel} links of alternating knots. A reverse parallel link of a knot consists of the two boundary components of an annulus $A$ embedded in $S^3$ with the said knot being one of the two components such that the two components are assigned opposite orientations. Let $K$ and $K^\p$ be the two components of a reverse parallel link induced by an annulus $A$. Following the convention that has been used in the literature (such as in \cite{Nutt,Rudolph}), we shall call the linking number $f$ between $K$ and $K^\p$ when they are assigned parallel orientations the {\em framing} of $K$. We note that a reverse parallel link of $K$ with framing $f$ is denoted by $K\ast_f A$ in \cite{Nutt} and by ${\rm{Bd}}_A(K,f)$ in \cite{Rudolph}. The framing is independent of the orientation of $K$ and the ambient isotopy class of $A$ in $S^3$ depends only on $K$ and the framing. Therefore, the reverse parallel links of $K$ are characterised by the framing $f$. Since our results (and proofs) only depend on the framing, not the actual annulus $A$, we shall introduce a new notation $\mathbb{K}_f$ for the reverse parallel link of $K$ with framing $f$. Keep in mind that the framing $f$ is the linking number of the two components of $\mathbb{K}_f$ with parallel orientations, hence the linking number of $\K_f$ itself is $-f$. 

\begin{figure}[!hbt]
\begin{center}
\includegraphics[scale=1.0]{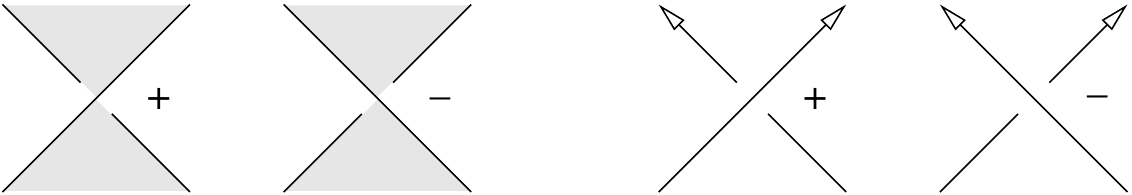}
\end{center}
\caption{Left: the crossing with respect to a checkerboard shading; Right: the crossing sign with respect to the orientation of the knot. 
\label{channel}}
\end{figure}

For a given knot diagram $D$ with a checkerboard shading, a crossing can be assigned a $+$ or a $-$ sign relative to this shading as shown on the left side of Figure \ref{channel}. This is not to be confused with the crossing sign with respect to the orientation of the knot which is used in the definition of the writhe of $D$ as shown on the right side of Figure \ref{channel}. 

\begin{figure}[!hbt]
\begin{center}
\includegraphics[scale=1.0]{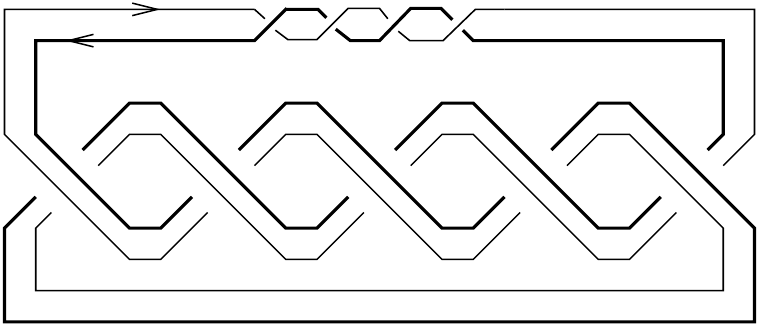}
\end{center}
\caption{The blackboard reverse parallel of the $(2,5)$ torus knot with two left-handed twists added. The framing of the resulting reverse parallel link (with the added twists) is thus $5+(-2)=3$. 
\label{blackboard}}
\end{figure}

Now, let $K$ be an alternating knot with a reduced diagram $D$. It is known that in such a case all crossings of $D$ are positive with respect to one checkerboard shading of $D$ and are all negative with respect to the other checkerboard shading of $D$. Furthermore, if we let $v_+(D)$ be the number of regions in the
shading with respect to which all crossings are positive, and $v_-(D)$ be the number of regions in the complementary shading with respect to which all crossings are negative, then $v_+(D)+v_-(D)-2=c(D)$  where $c(D)$ is the number of crossings in $D$ \cite{Kauﬀman1987}. From $D$ we can also obtain its so-called 
blackboard reverse parallel annulus (framing) which provides a good reference for other choices of annuli (framings) as the other choices come from this one by adding either right-handed or left-handed twists. If the writhe of $D$ is $w(D)$, then the framing of the blackboard reverse parallel is also $w(D)$. If $k$ right-handed (left-handed) twists are added between the two components, then the resulting reverse parallel has framing $w(D)+k$ ($w(D)-k$). See Figure \ref{blackboard} for an illustration.

Our main result in this paper is the following theorem.

\begin{theorem}\label{main}
Let $K$ be an alternating knot and $D$ a reduced diagram of $K$. Let $c(D)$, $w(D)$, $v_+(D)$ and $v_-(D)$ be as defined above, then the braid index of $\mathbb{K}_f$, denoted by $\b(\mathbb{K}_f)$, is given by the following formula
\begin{equation}\label{MainE}
\b(\mathbb{K}_f)=\left\{
\begin{array}{ll}
c(D)+2+a(D)-f, &\ {\rm if}\ f < a(D),\\
c(D)+2, &\ {\rm if}\ a(D)\le f \le b(D),\\
c(D)+2-b(D)+f, &\ {\rm if}\ f > b(D),
\\
\end{array}
\right.
\end{equation}
where $a(D)=-v_-(D)+w(D)$ and $b(D)=v_+(D)+w(D)$.
\end{theorem}

We can summarise theorem \ref{main} pictorially in terms of the blackboard reverse parallel of $D$.
\begin{itemize}\item The blackboard reverse parallel has braid index $c(D)+2$.
\item The braid index remains  $c(D)+2$ after adding up to $v_+(D)$ right-hand twists, or up to $v_-(D)$ left-hand twists.
\item Each further right or left-hand twist increases the braid index by $1$.
\end{itemize}
So for example, since $v_-(D)=2$ and $v_+(D)=5$ for the $(2,5)$ torus knot, the braid index for the reverse parallel shown in figure \ref{blackboard}  is $c(D)+2=7$.  Adding one further left-hand twist would increase the braid index to $8$, while we would still have braid index $7$ after adding anything up to $5$ right-hand twists to the blackboard parallel.

We shall establish (\ref{MainE}) by proving that the right side expression is both a lower bound and an upper bound for the $\b(\mathbb{K}_f)$. The lower bound is obtained by the Morton-Franks-Williams inequality while the upper bound is established by direct construction.

\section{The lower bound}\label{s2}

In this section, we shall prove the following theorem.
\begin{theorem}\label{lowerbound}
Let $\mathbb{K}_f$ be the reverse parallel link of an alternating knot $K$ with framing $f$ and $D$  a reduced diagram of $K$.  Then\begin{equation}\label{lowerbound_formula}
\b(\mathbb{K}_f)\ge \left\{
\begin{array}{ll}
c(D)+2+a(D)-f, &\ {\rm if}\ f < a(D),\\
c(D)+2, &\ {\rm if}\ a(D)\le f \le b(D),\\
c(D)+2-b(D)+f, &\ {\rm if}\ f > b(D).\\
\end{array}
\right.
\end{equation}
where $a(D)=-v_-(D)+w(D)$ and $b(D)=v_+(D)+w(D)$.
\end{theorem}

\subsection{The Homfly and Kauffman polynomials}
Before proving this theorem we note some properties of the Homfly and Kauffman polynomials of a link $L$.

 The \emph{Homfly polynomial} $P_L(v,z)\in \mathbb{Z}[v^{\pm1},z^{\pm1}]$ of an oriented link $L$ is  determined by the skein relations
 \[v^{-1}P_{L^+}-vP_{L^-}=z P_{L^0}\] where $L^{\pm}, L^0$ differ only near one crossing as shown and    takes the value $1$ on the unknot.
 \[L^+ =\Xor, \quad L^-=\Yor, \quad L^0=\Ior.\]
 
 The \emph{Kauffman polynomial} $F_L(a,z)\in \mathbb{Z}[a^{\pm1},z^{\pm1}]$ for an unoriented link $L$  is defined in \cite{Kauffman1990}. Again it takes the value $1$ on the unknot.
 
 When an extra distant unknotted component $O$ is adjoined to the link $L$, to make $L\sqcup O$, each polynomial changes in the following simple way.
 \beqn
 P_{L\sqcup O}(v,z)&=&\frac{v^{-1}-v}{z}P_L(v,z)\\
 F_{L\sqcup O}(a,z)&=& \left(\frac{a+a^{-1}}{z}-1\right)F_L(a,z).
 \eeqn
 
Define the \emph{extended Homfly polynomial} $EP$ by
\begin{equation}\label{EP_Def}
EP_L(v,z)=\frac{v^{-1}-v}{z}P_L(v,z)=P_{L\sqcup O}(v,z)
\end{equation}  
and the \emph{extended Kauffman polynomial }$EF$ by 
\begin{equation}\label{EF_Def}
EF_L(a,z)=\left(\frac{a+a^{-1}}{z}-1\right)F_L(a,z)=F_{L\sqcup O}(a,z).
\end{equation}  

\begin{remark} This extended normalisation is often used in the context of quantum invariants, where it allows for more natural specialisations of the knot polynomials. It is also more useful in that context to use the Dubrovnik variant of the Kauffman polynomial in place of $F$. 

 By plugging in $L=\phi$ on both sides of (\ref{EP_Def}) and (\ref{EF_Def}) the extended polynomials can be thought of as taking the value $1$ on the empty link $\phi$.

 \end{remark}

 \subsection{Bounds from the Homfly and Kauffman polynomials}
 
 The Morton-Franks-Williams inequality \cite{Franks,Morton1986} gives a lower bound for the braid index $\b(L)$ of the link $L$ in terms of the $v$-spread of the Homfly polynomial $P_L(v,z)$ or its extended version. Explicitly
 \begin{equation}\label{MFW}
\b(L)\ge 1+\frac{1}{2}{\rm spr}_vP_{L}(v,z)=\frac{1}{2}{\rm spr}_vEP_{L}(v,z).
\end{equation}

The $a$-spread of the Kauffman polynomial is shown by Morton and Beltrami \cite{Morton1998} to give a bound for the arc index $\alpha(L)$.
Explicitly this is \[{\rm spr}_aF_L(a, z)\le \alpha(L)-2.\]
Bae and Park \cite{Bae}  showed that the arc index  $\alpha(L)$ is bounded above by $c(L)-2$, that is, $\alpha(L)\le c(L)+2$. Combining these results shows that
\begin{equation} {\rm spr}_aF_L(a, z)\le c(L). \label{arcindex}
\end{equation}

\subsection{A congruence result} 

Rudolph  \cite{Rudolph} relates the Kauffman polynomial of a link $L$ with the Homfly polynomial of the reverse parallels  of $L$. 

\begin{notation}
For Laurent polynomials $A =\sum a_{i,j}v^iz^j, B =\sum b_{i,j}v^iz^j \in  \mathbb{Z}[v^{\pm1},z^{\pm1}]$  we write  $A\cong_{\mathbb Z_2} B$ when  $a_{i,j}\cong b_{i,j} \bmod 2$ for all $i,j$.
\end{notation}

 In the case of a knot $K$ Rudolph's theorem for the reverse parallel  $\mathbb{K}_f$  can then be stated very cleanly in terms of the extended polynomials.

\begin{theorem}\cite[Congruence Theorem]{Rudolph}\label{Cong_Theorem}
\[EP_{\mathbb{K}_f}(v,z) -1 \cong_{\mathbb{Z}_2} v^{-2f}EF_K(v^{-2},z^2)\]
\end{theorem}

\subsection{Alternating knots}
We can apply these bounds to the case of alternating knots, starting from observations of Cromwell \cite{Cromwell} about their Kauffman polynomial.

For any knot $K$ with a diagram $D$,  write the Kauffman polynomial $F_K(a,z)$ of $K$ as
\begin{equation}\label{e2.3}
F_K(a, z) = a^{-w(D)}\sum_{i,j}a_{i,j}a^iz^j.
\end{equation}
In this form the coefficients $a_{i,j}$ are only non-zero in the range $|i| + j\le c(D)$. 

Cromwell extends work of Thistlethwaite    \cite{Thistlethwaite} to identify two non-zero coefficients $a_{i,j}$ which realise the maximum possible $a$-spread $c(D)$ for $F_K(a,z)$ in the case of an alternating knot $K$ with reduced diagram $D$.

\begin{theorem}\cite{Cromwell}\label{TCromwell} 
Let $K$ be an alternating knot and $D$  a reduced diagram of $K$. Then
$a_{i,j} =1$ in the two cases $i=1-v_+(D), j=c(D)+i$ and  $i=v_-(D)-1,j=c(D)-i$. 
\end{theorem}

It follows that ${\rm spr}_aF_K(a, z)\ge v_-(D)-1-(1-v_+(D))= c(D)$.
\begin{corollary}\label{Cor2.4}
Hence ${\rm spr}_aF_K(a, z)= c(D)$, by (\ref{arcindex}). Then $a_{i,j}=0$ in (\ref{e2.3}) unless $1-v_+(D)\le i\le v_-(D)-1$.
\end{corollary}

Now set
\begin{equation}\label{B}
B_D(a,z)=a^{w(D)}EF_K(a,z)=\left(\frac{a+a^{-1}}{z}-1\right)\sum_{i,j}a_{i,j}a^iz^j.\end{equation}
Then ${\rm spr}_aB_D(a,z)={\rm spr}_aF_K(a, z)+2=c(D)+2$. Furthermore, if we write 
\begin{equation}\label{B2}
B_D(a,z)=\sum_{i,j}b_{i,j}a^iz^j,
\end{equation} then $b_{i,j} =0 $ unless $-v_+(D)\le i\le v_-(D)$ by  corollary \ref{Cor2.4}. 

The two critical  monomials  $a^{-v_+(D)}z^{c(D)-v_+(D)}$ and  $a^{v_-(D)}z^{c(D)-v_-(D)}$ in $B_D(a,z)$, which correspond to $i=-v_+(D)$ 
and $i=v_-(D)$ respectively, both have coefficient $b_{i,j}=1$, by theorem \ref{TCromwell}. 
We will use these critical monomials in finding a 
 lower bound for the $v$-spread of the extended Homfly polynomial of the  reverse parallels of $D$.
 
 Theorem \ref{e2.5} below gives a simple formula to calculate the extended Homfly polynomial of $\mathbb{K}_{k+f}$ in terms of the polynomial of $\mathbb{K}_{k}$. 
 
 
 
\begin{theorem}\label{e2.5}
For any $f$ and $k$ we have
\[v^{2f}(EP_{\mathbb{K}_{k+f}}(v,z)-1)=EP_{\mathbb{K}_{k}}(v,z)-1.\]
\end{theorem}

\begin{proof}
While this is in effect shown by Rudolph in his proposition 2(5) it is easy to give a direct skein theory proof.
It is enough to prove in the case $f=1$.
Now $\mathbb{K}_{k+1}$ is given from $\mathbb{K}_{k}$ by adding one extra twist in the annulus, as shown.  
\[\mathbb{K}_{k}=\labellist\small
\pinlabel {$K$} at 65 90
\endlabellist \kparallelstraight, \quad \mathbb{K}_{k+1} = \labellist\small
\pinlabel {$K$} at 65 90
\endlabellist \kparallel.\]
With the reverse parallel orientation on the strings apply the Homfly skein relation at one of the crossings in the diagram for  $\mathbb{K}_{k+1}$. Since this is a negative crossing $\mathbb{K}_{k+1}$ plays the role of $L^-$. Switching the crossing gives \[L^+ =\labellist\small
\pinlabel {$K$} at 65 90
\endlabellist \kparallelswitch = \labellist\small
\pinlabel {$K$} at 65 90
\endlabellist \kparallelstraight \quad=\ \mathbb{K}_{k}\] while the smoothed diagram \[L^0=\labellist\small
\pinlabel {$K$} at 65 90
\endlabellist \kparallelsmooth\] is simply an unknotted curve.

The skein relation, in the form
\[EP_{L^+}=vz EP_{L^0}+v^2EP_{L^-},\] then gives
\beqn EP_{\mathbb{K}_{k}}&=&vz\frac{v^{-1}-v}{z}+v^2EP_{\mathbb{K}_{k+1}}\\
&=&1-v^2 +v^2EP_{\mathbb{K}_{k+1}}
\eeqn
Thus \[v^{2}(EP_{\mathbb{K}_{k+1}}-1)=EP_{\mathbb{K}_{k}}-1.\]
\end{proof}

We can now specify a lower bound for the $v$-spread of the extended Homfly polynomial of the parallels $\mathbb{K}_{w(D)+f}$ as $f$ varies.
\begin{theorem}\label{spread}
Let $K$ be an alternating knot with reduced diagram $D$.
The framed reverse parallel $\mathbb{K}_{w(D)+f}$ has the following lower bound for the $v$-spread of its extended Homfly polynomial.
\[{\rm spr}_vEP_{\mathbb{K}_{w(D)+f}}(v,z)\ge 
\left\{
\begin{array}{ll}
2(v_+(D)-f), &\ {\rm if}\ f < -v_-(D),\\
2(v_+(D)+v_-(D)), &\ {\rm if}\ -v_-(D)\le f \le v_+(D),\\
2(f+v_-(D)), &\ {\rm if}\ f > v_+(D).
\\
\end{array}
\right.
\]
\end{theorem}

\begin{proof}
Where $K$ is an alternating knot with reduced diagram $D$ theorem \ref{Cong_Theorem} shows that
\begin{equation}\label{econgruence} 
B_D(v^{-2},z^2)=v^{-2w(D)}EF_K(v^{-2},z^2)\cong_{\mathbb{Z}_2} EP_{\mathbb{K}_{w(D)}}(v,z) -1. \end{equation}

Now in $B_D(v^{-2},z^2)=\sum b_{i,j}v^{-2i}z^{2j}$ there are two critical monomials $v^{-2i}z^{2j}$, one with $i=-v_+(D), j=c(D)-v_+(D)$ and the other with $i=v_-(D), j=c(D)-v_-(D)$, where  $b_{i,j}=1$. By equation (\ref{econgruence}) there are two corresponding critical monomials $v^{-2i}z^{2j}$ in $EP_{\mathbb{K}_{w(D)}}(v,z) -1$ whose coefficients are congruent to $b_{i,j}$, and hence are odd. One term has $v$-degree $-2v_-(D)$ and the other has $v$-degree $2v_+(D)$.

By theorem \ref{e2.5} we have
\[v^{2f} EP_{\mathbb{K}_{w(D)+f}}(v,z)=(EP_{\mathbb{K}_{w(D)}}(v,z) -1) +v^{2f}.\]

The $v$-spread of $EP_{\mathbb{K}_{w(D)+f}}(v,z)$ is the same as the $v$-spread of $(EP_{\mathbb{K}_{w(D)}}(v,z) -1) +v^{2f}$.
In this Laurent polynomial 
consider the appearance of the two critical monomials along with the monomial $v^{2f}$. Unless one of the two critical monomials $v^{2v_+(D)}z^{2c(D)-2v_+(D)}$ and $v^{-2v_-(D)}z^{2c(D)-2v_-(D)}$ in $B_D(v^{-2},z^2)$ is $v^{2f}$ they will each still have odd coefficients, and the $v$-spread will be at least $2(v_+(D)+v_-(D))$.

If $f<-v_-(D)$ or $f>v_+(D)$ the monomial $v^{2f}$ has even coefficient in  $EP_{\mathbb{K}_{w(D)}}(v,z) -1$ since it has coefficient $0$ in $B_D(v^{-2},z^2)$. In this range of $f$ it then has non-zero coefficient in
$(EP_{\mathbb{K}_{w(D)}}(v,z) -1) +v^{2f}$. This gives the lower bound $2(v_+(D)-f)$ when $f<-v_-(D)$, and  $2(v_-(D)+f)$ when $f>v^+(D)$ for ${\rm spr}_vEP_{\mathbb{K}_{w(D)+f}}(v,z)$. 

To complete the proof of theorem \ref{spread} it remains to deal with the cases where $v^{2f}$ is one of the two critical monomials $v^{2v_+(D)}z^{2c(D)-2v_+(D)}$ and $v^{-2v_-(D)}z^{2c(D)-2v_-(D)}$ in $B_D(v^{-2},z^2)$. In the first case this means that $f=v_+(D)$ and $0=c(D)-v_+(D)$. Then $f=c(D)=v_+(D)=n$ and $D$ is the reduced diagram of the $(2,n)$ torus knot. In the other case $-f=c(D)=v_-(D)=n$ hence $D$ is the reduced diagram of the $(2,-n)$ torus knot.

In the case that $D$ is the $(2,n)$ torus knot, we need to show that the coefficient of $v^{2n}$ in $(EP_{\mathbb{K}_{w(D)}}(v,z) -1) +v^{2n}$ is non-zero. In theorem \ref{torus} we show that this coefficient is $2$, by showing that $v^{2n}$ has coefficient $1$ in $EP_{\mathbb{K}_{w(D)}}(v,z)$, where $\mathbb{K}_{w(D)}$ is the blackboard reverse parallel of $D$.

The  case of the $(2,-n)$ torus knot follows directly by considering the polynomial of the mirror image and this completes the proof of theorem \ref{spread}. \end{proof}

The detailed calculation for the special case of the $(2,n)$ torus knot will now be shown. 

\begin{theorem} \label{torus} The blackboard reverse parallel $\mathbb{K}_n$ of the $(2,n)$ torus knot $K$ satisfies
 \[EP_{\mathbb{K}_n}(v,z)=v^{2n} + \sum_{i<2n, j} a_{i,j} v^i z^j.\]
 \end{theorem}

\begin{proof} We can draw a diagram of $\mathbb{K}_n$ as the closure of a $4$-strand tangle with two upward and two downward strings, as shown. 
\[ \torusparallel\]
 It is  more convenient to place the upward pair of strings at the left, at the top and bottom, and write $\mathbb{K}_n$ as the closure of the tangle $T^n$, where \[T= \twotwotangle\]
 
 We use the skein relations in the form
  \[v^{-1}\Xor-v\Yor \ =\ z \Ior\] to write the closure of $T^n$ as a linear combination of the closures of simpler tangles.

  \begin{notation}
  We will say that the $4$-strand tangle $U$ \emph{evaluates to} the extended Homfly polynomial of its closure, which we write as $\mathrm{ev}(U)\in  \mathbb{Z}[v^{\pm1},z^{\pm1}]$.
  \end{notation}
   \begin{remark}
   Evaluation is linear on tangles, and respects the skein relations. It is a sort of trace function in that $\ev(AB)=\ev(BA)$.
   \end{remark}
  
  Our first step is to expand $T$ as a combination of the tangles \[\sigma_1= \twotwosigmaone\ , \ \sigma_3 =\twotwosigmathree\ ,\  h=\twotwoh\ , \ H=\twotwohh\ , \] and their products when placed one above the other.
  \begin{remark} By using the skein relations we are in effect working in a version of the mixed Hecke algebra $H_{2,2}(v,z)$ spanned by tangles with two upward and two downward strings \cite{Morton2006}.
  \end{remark}
 
 The crossing circled here   in \[T \quad=\quad \twotwotanglecircle\]  is a negative crossing so we can use the skein relation at this crossing in the form \[\Yor= v^{-2}\Xor-v^{-1}z\Ior .\]
  
  Then we have 
\beqn T=\twotwotangle&=&v^{-2}\ \Tswitch \quad -\quad v^{-1}z\ \Tsmooth\\[1mm]
&=& v^{-2}\sigma_1\sigma_3 -v^{-1}z\sigma_1\sigma_3h\\
&=&C+C\tau,
\eeqn
 when for convenience we set $C=c_1 c_3=(v^{-1}\sigma_1)(v^{-1}\sigma_3)$ and $\tau=(-zv)h$.
 
 Then $T^n =(C+C\tau)^n$. Now $C$ and $ \tau$ do not commute so we write
 \begin{equation}\label{Tnexpansion}T^n=C^n +(C\tau)^n +\sum_{0<k<n} C^{r_1}\tau C^{r_2}\tau \cdots C^{r_k}\tau C^r,\end{equation}
 where $r_i\ge 1,r\ge 0$ and $r+\sum r_i =n$.

 We can estimate the contribution of these terms to the evaluation of $T^n$.
 
 \begin{itemize}
\item The evaluation of $C^n$ only contributes terms up to $v$-degree $4$.

\item The terms in the large sum with weight $k$ in $\tau$ evaluate to terms of $v$-degree at most $2k$. Without changing the evaluation we can assume that $r=0$, since we can cycle $C^r$ from the end to the beginning of the product, and amalgamate it with $C^{r_1}$. The contribution of these terms with $k<n$ to the evaluation of $T^n$ is shown in proposition \ref{kterms} to have degree no more than $2n-2$ in $v$.

\item The most important contribution comes from the evaluation  of $(C\tau)^n$, which gives $v^{2n}$, and no other terms with $v$-degree $2n$ or larger, as stated in proposition \ref{prophighest}.  
\end{itemize}

Before making detailed calculations we note
some useful properties which can be quickly checked diagrammatically.
\begin{itemize}
\item $\sigma_1 H=\sigma_3H,\ H\sigma_1=H\sigma_3$
\item $H=h\sigma_1\sigma_3^{-1}h $
\item $\quad\Idor  \unknot  \ =\  \delta\  \Idor $, where $\delta=\frac{v^{-1}-v}{z}$
\item $h^2=\delta h$
\item $h\sigma_1h=h$
\end{itemize}

Here are some consequences for our use of $c_1=v^{-1}\sigma_1,c_3 =v^{-1}\sigma_3, C=c_1c_3,\tau =(-zv)h$, which  follow algebraically.
\begin{itemize}
\item $c_1=c_1^{-1}+z, \quad c_3=c_3^{-1}+z$, (skein relation)
\item $\tau c_1 c_3^{-1}\tau=(-zv)^2hc_1c_3^{-1}h=(zv)^2H$
\item $\tau^2=(-zv)\delta \tau=(v^2-1)\tau$
\item $\tau c_1\tau=v^{-1}(-zv)^2 h\sigma_1 h= -z\tau$
\item $\tau C\tau=\tau(c_1c_3^{-1}+zc_1)\tau=(zv)^2H-z^2\tau$
\end{itemize}

\medskip
\begin{proposition} \label{prophighest}
The extended polynomial of the closure of $(C\tau)^n$ is $v^{2n}+\text{ lower terms in }v$, for $n>1$, and $1-v^{-2}$ when $n=1$.
\end{proposition}

\begin{proof}[Proof of proposition]
When $n=1$ we have $C\tau =(-zv^{-1})\sigma_1\sigma_3 h$. Now $\sigma_1\sigma_3h$ closes to a single unknotted curve, so $C\tau$ evaluates to $-zv^{-1}\delta=1-v^{-2}$.

For $n>1$ write \beqn (C\tau)^n &=&C(\tau C \tau) (C\tau)^{n-2}\\
&=& (zv)^2 C H (C\tau)^{n-2} -z^2 C\tau (C\tau)^{n-2}
\eeqn
The evaluation of the second term has  $v$-degree at most $2n-2$, by induction on $n$, so any monomials of larger $v$-degree  must come from the first term.

Now $Hh=\delta H$ and $H\sigma_1 h=H$
We can then write \beqn HC\tau &=& H(c_1c_3^{-1}+zc_1)\tau \\
&=&H\tau+z Hc_1\tau\\
&=&(-zv)(\delta+zv^{-1})H\\
&=&(v^2-1-z^2)H.
\eeqn
So the first term expands to \[(zv)^2 C H (C\tau)^{n-2} =(zv)^2(v^2-1-z^2)^{n-2}CH.\]
 Now $CH=c_1c_3^{-1}H +zc_1H=H+zv^{-1}\sigma_1H$. The closure of $H$ is two disjoint unknotted curves, and $\sigma_1H$ closes to one unknotted curve, evaluating to $\delta^2$ and $\delta$ respectively.
The first term then evaluates to \[(v^2-1-z^2)^{n-2}(\delta^2(-zv)^2-z^2(-zv\delta ))=(v^2-1-z^2)^{n-1}(v^2-1).\] This contributes a single term $v^{2n}$ and no further terms of $v$-degree larger than $2n-2$.\end{proof}

We now show that the remaining terms in the expansion of $T^n$ in (\ref{Tnexpansion})
contribute terms  in $v$ of degree $\le 2n-2$, when $n\ge 3$.

The skein relation, in the form $c_1^2=1+zc_1$
allows us to write $c_1^r$ recursively as a linear combination of $c_1$ and the identity tangle 
\[c_1^r=a_r(z)+b_r(z)c_1\] with coefficients which are polynomials in $z$ only. 

We can then expand $C^r$ as a linear combination of $C, c_1,c_3$ and the identity tangle, with coefficients in $\mathbb{Z}[z]$. Explicitly
\[C^r=(a_r+b_rc_1)(a_r+b_rc_3)=a_r^2+a_rb_r(c_1+c_3)+b_r^2C.\]

Now $\sigma_1\sigma_3$ closes to two unknotted curves, evaluating to $\delta^2$, $\sigma_1$ and $\sigma_3$ close to three unknots and the identity tangle closes to four unknots.

The term $C^n$ in the expansion of $T^n$ then contributes $(a_n\delta^2+b_nv^{-1}\delta)^2$ to the evaluation. This provides terms of degree at most $4 $ in $v$.

To complete our proof of theorem \ref{torus} we show that the evaluation of the remaining terms in (\ref{Tnexpansion})
 has $v$-degree at most $2n-2$.

 This follows from
\begin{proposition}\label{kterms}
The evaluation  of \[ C^{r_1}\tau\cdots C^{r_i}\tau\cdots C^{r_k}\tau\] with $r_i\ge 1$  has terms of degree at most $2k$ in $v$.
 \end{proposition}
 \begin{proof}

 By induction on the number of exponents $r_i$ for which $r_i>1$.
 
  When $r_i=1$ for all $i$ this follows from proposition \ref{prophighest}.
 
 Otherwise  we can cycle the terms in the product without changing its evaluation, and arrange that $r_k=r>1$. 
 Then \beqn \tau C^r \tau&=&a_r^2 \tau^2+a_rb_r\tau(c_1+c_3)\tau +b_r^2\tau C\tau\\
 &=&a_r^2(v^2-1)\tau -2z a_rb_r\tau +b_r^2\tau C\tau.
 \eeqn
 Then \beqn C^{r_1}\tau\cdots  C^{r_k}\tau&=&\left(a_r^2(v^2-1)-2z a_rb_r\right)
 C^{r_1}\tau\cdots C^{r_{k-1}}\tau\\
 && +b_r^2\ C^{r_1}\tau\cdots  C^{r_{k-1}}\tau C\tau\eeqn
 
These  expressions both have one less term $C^{r_i}$ for which $r_i>1$, so by our induction hypothesis the evaluation of $C^{r_1}\tau\cdots  C^{r_{k-1}}\tau C\tau$ has terms of degree at most $2k$ in $v$ while $C^{r_1}\tau\cdots C^{r_{k-1}}\tau $ has terms of degree at most $2k-2$. With the coefficient  $a_r^2(v^2-1)-2z a_rb_r$ adding $2$  in this case all terms in the final evaluation have degree at most $2k$ in $v$.  This establishes the proposition.
 \end{proof}
 Now all the terms in (\ref{Tnexpansion}) have been dealt with, and 
theorem \ref{torus} for the evaluation of the reverse  blackboard parallel of the $(2,n)$ torus knot then follows.
\end{proof}
The proof of theorem \ref{spread} is then complete.
We can now prove theorem \ref{lowerbound} which was the goal of this section.

\begin{proof}[Proof of theorem \ref{lowerbound}]
Using the Morton-Franks-Williams bound (\ref{MFW}) in theorem \ref{spread} immediately gives the  lower bound for the braid index of $\mathbb{K}_{w(D)+f}$ as
\[\b(\mathbb{K}_{w(D)+f})\ge 
\left\{
\begin{array}{ll}
v_+(D)-f, &\ {\rm if}\ f < -v_-(D),\\
v_+(D)+v_-(D), &\ {\rm if}\ -v_-(D)\le f \le v_+(D),\\
f+v_-(D), &\ {\rm if}\ f > v_+(D).
\\
\end{array}
\right.
\]
Replacing  $f$ by $f-w(D)$ then gives
\[\b(\mathbb{K}_{f})\ge 
\left\{
\begin{array}{ll}
v_+(D)-f+w(D), &\ {\rm if}\ f-w(D) < -v_-(D),\\
v_+(D)+v_-(D), &\ {\rm if}\ -v_-(D)\le f -w(D)\le v_+(D),\\
f-w(D)+v_-(D), &\ {\rm if}\ f -w(D)> v_+(D).
\\
\end{array}
\right.
\]
Now $v_+(D)+v_-(D)=c(D)+2$, so after setting $a(D)=w(D)-v_-(D)$ and $b(D)=w(D)+v_+(D)$ this lower bound becomes
 \[\b(\mathbb{K}_{f})\ge 
 \left\{
\begin{array}{ll}
c(D)+2+a(D)-f, &\ {\rm if}\ f < a(D),\\
c(D)+2, &\ {\rm if}\ a(D)\le f \le b(D),\\
c(D)+2-b(D)+f, &\ {\rm if}\ f > b(D).\\
\end{array}
\right .
\]
which is  the formula (\ref{lowerbound_formula}) claimed in  theorem \ref{lowerbound}
\end{proof}

\section{The upper bound}\label{s3}

In this section, we shall prove the following theorem, which provides us the desired upper bound for the braid index of $\mathbb{K}_f$.

\begin{theorem}\label{upperbound}
If $\mathbb{K}_f$ is a reverse parallel link of an alternating knot $K$ with framing $f$ and $D$ is a reduced diagram of $K$, then we have
\begin{equation}\label{upperbound_formula}
\b(\mathbb{K}_f)\le \left\{
\begin{array}{ll}
c(D)+2+a(D)-f, &\ {\rm if}\ f < a(D),\\
c(D)+2, &\ {\rm if}\ a(D)\le f \le b(D),\\
c(D)+2-b(D)+f, &\ {\rm if}\ f > b(D).\\
\end{array}
\right.
\end{equation}
where $a(D)=-v_-(D)+w(D)$ and $b(D)=v_+(D)+w(D)$.
\end{theorem}

\begin{proof} It suffices to show that braid presentations of $\mathbb{K}_f$ can be constructed with the number of strings given in the theorem. The construction depends on using an arc presentation for $K$. From
an arc presentation for $K$ with $\alpha$ arcs, Nutt \cite{Nutt} gives an upper bound for $\b(\mathbb{K}_f)$ in the following form 
\begin{equation}\label{Nutt_formula}
\b(\mathbb{K}_f)\le \left\{
\begin{array}{ll}
\alpha+a^\prime-f, &\ {\rm if}\ f < a^\prime,\\
\alpha, &\ {\rm if}\ a^\prime\le f \le b^\prime,\\
\alpha+f-b^\prime, &\ {\rm if}\ f > b^\prime,\\
\end{array}
\right.
\end{equation}
where $a^\prime$, $b^\prime$ are some integers. We need to refine (\ref{Nutt_formula}) for our purpose. First we make the following observation. In his argument leading to the proof of \cite[Theorem 3.1]{Nutt}, Nutt constructed a closed braid template for $\mathbb{K}_f$. In this template there are $\alpha-1$ handles. Placing a single crossing in each handle will result in a closed braid representation of $\mathbb{K}_f$ for some framing $f$ with $\alpha$ strings. By starting with the template in which a positive crossing is placed in every handle, then switching the sign of a positive crossing, one at a time, until all crossings are negative, we then obtain a sequence of $\alpha$ closed braid representations $\beta_1$, $\beta_2$, ..., $\beta_\alpha$ (each with $\alpha$ strings) of reverse parallel links 
$\mathbb{K}_{f_1}$, $\mathbb{K}_{f_2}$, ..., $\mathbb{K}_{f_\alpha}$ for some framings $f_1$, $f_2$, ..., $f_\alpha$ with the property $f_{j+1}=f_j-1$, $j=1, ..., \alpha-1$. Combining this observation with the fact that $K$ has an arc representation with $c(D)+2$ arcs \cite{Bae}, (\ref{Nutt_formula}) can then be improved to the following inequality.
\begin{equation}\label{Nutt_formula2}
\b(\mathbb{K}_f)\le \left\{
\begin{array}{ll}
c(D)+2+a^\prime-f, &\ {\rm if}\ f < a^\prime,\\
c(D)+2, &\ {\rm if}\ a^\prime\le f \le b^\prime,\\
c(D)+2-b^\prime+f, &\ {\rm if}\ f > b^\prime,\\
\end{array}
\right.
\end{equation}
where $b^\prime=c(D)+2+a^\prime$. 

We now claim that we must have $a^\prime=a(D)$ and $b^\prime=b(D)$. Say $a^\prime>a(D)$. Then 
$$
b^\prime=a^\prime+c(D)+2>a(D)+c(D)+2=b(D).
$$ 
Thus for $f=b^\prime>b(D)$, we have 
$$
\b(\mathbb{K}_f)\ge c(D)+2+b^\prime-b(D)>c(D)+2
$$ 
by (\ref{upperbound_formula}) and 
$\b(\mathbb{K}_f)\le c(D)+2$ by (\ref{Nutt_formula2}), which is a contradiction. On the other hand, if $a^\prime<a(D)$, then for  $f=a^\prime<a(D)$, we have
$$
\b(\mathbb{K}_f)\ge c(D)+2+a(D)-a^\prime>c(D)+2
$$ 
by (\ref{upperbound_formula}) and 
$\b(\mathbb{K}_f)\le c(D)+2$ by (\ref{Nutt_formula2}), which is again a contradiction. This shows that $a^\prime=a(D)$. Similarly, we can prove that $b^\prime=b(D)$. This concludes the proof of Theorem \ref{upperbound}.
\end{proof}

We end our paper with the following remarks.

  \begin{remark}
 In the case that one desires using the linking number $\ell$ of $\K_f$ in the formulation of $\b(\K_f)$ instead of the framing $f$, then the formulation can be easily obtained by substituting $f$ by $-\ell$ in (\ref{MainE}). Specifically, (\ref{MainE}) becomes
 \begin{equation}\label{MainE2}
\b(\mathbb{K}_f)=\left\{
\begin{array}{ll}
c(D)+2+a^\p(D)-\ell, &\ {\rm if}\ \ell < a^\p(D),\\
c(D)+2, &\ {\rm if}\ a^\p(D)\le \ell \le b^\p(D),\\
c(D)+2-b^\p(D)+\ell, &\ {\rm if}\ \ell > b^\p(D),
\\
\end{array}
\right.
\end{equation}
where $a^\p(D)=-v_+(D)-w(D)$, $b^\p(D)=v_-(D)-w(D)$ and $\ell=-f$ is the linking number of $\K_f$. The corresponding formulation of (\ref{lowerbound_formula}) matches the one given in \cite{Diao2021_1}. We need to point out that the lower bound formula derived in \cite{Diao2021_1} uses a graph theoretic approach on the Seifert graphs of $D$ and $\K_f$ constructed from the blackboard reverse parallel of $D$. However that approach only works for the special alternating knots, namely those alternating knots which admit a reduced alternating diagram in which the crossings are either all positive or all negative.
  \end{remark}
  
  \begin{remark}
  The general question of finding the braid index for a satellite of a knot $K$ with some form of reverse string pattern has been considered by Birman and Menasco \cite{Birman}. Our reverse parallels, along with Whitehead doubles, are the simplest such satellites. Nutt \cite{Nutt} draws on \cite{Birman} to give lower bounds for the braid index in terms of the arc index of $K$, as well as the upper bounds which we have used. Coupled with the later work of Bae and Park \cite{Bae} this would provide our result without the use of Rudolph's congruence.

Some descriptions given by \cite{Birman} were later found to be incomplete, with Ka Yi Ng \cite{Ng} providing details of the missing cases. Nutt's lower bound argument  needs  the analysis in \cite{Birman} which shows  that the arc index of $K$ is a lower bound for the braid index of any reverse string satellite of $K$.
  We have not been able to confirm how well the arc index analysis in the original paper extends to Ng's extra cases.

  \end{remark}

  \begin{remark}
  Theorem \ref{main} allows us to settle a long standing conjecture regarding the ropelength of an alternating knot. The conjecture states that the ropelength of an alternating knot $K$ is at least proportional to its crossing number. This statement is now a consequence of \cite[Theorem 3.1]{Diao2020} and the fact that the ropelength of $K$ is bounded below by a (fixed) constant multiple of the ropelength of $\K_f$ for some $f$. What remains open is the more general case of an alternating link with two or more components.
  \end{remark}

\end{document}